# Optimization Method for Interval Portfolio Selection Based on Satisfaction Index of Interval inequality Relation

**Yunchol Jong**[a]

[a]Center of Natural Science, University of Sciences, Pyongyang, DPR Korea

E-mail: yuncholjong@yahoo.com

**Abstract:** In this paper we consider an interval portfolio selection problem with uncertain returns and introduce an inclusive concept of satisfaction index for interval inequality relation. Based on the satisfaction index, we propose an approach to reduce the interval programming problem with uncertain objective and constraints into a standard linear programming problem with two parameters. We showed by simulation experiment that our method is capable of helping investors to find efficient portfolios according to their preference.

## 1. Introduction

The theory of optimal portfolios selection was developed by Harry Markowitz in the 1950's. His work [1] formalized the diversification principle in portfolio selection and earned him the 1990 Nobel Prize for Economics. Consider an investor who has a certain amount of money to be invested in a number of different securities (stocks, bonds, etc.) with uncertain returns. The portfolio vector $x$ must satisfy the first constraint $\sum_{i=1}^{n} x_i = 1, x_i \geq 0, i = 1,...,n$ and there may or may not be additional feasibility constraints. A feasible portfolio $x$ is called efficient if it has the maximal expected return among all portfolios with the same variance, or alternatively, if it has the minimum risk (variance) among all portfolios that have at least a certain expected return. The collection of efficient portfolios forms the efficient frontier of the portfolio universe. Markowitz' portfolio optimization problem, also called the mean-variance optimization (MVO) problem, can be formulated in three different but equivalent ways. One formulation (MVO1) results in the problem of finding a minimum variance portfolio of the securities 1 to $n$ that yields at least a target value R of expected return. Then the second constraint indicates that the expected return is no less than the target value and, as we discussed above, the objective function corresponds to the total risk of the portfolio. Nonnegativity constraints on $x_i, i = 1,...,n$ are introduced to rule out short sales (selling a security that you do not have). As an alternative to the problem MVO1, we may consider the problem (MVO2) that is to maximize the expected return of a portfolio while limiting the variance of its return.

Traditionally, it has been assumed that the distribution functions of possibility returns are known while solving portfolio selection models. However, new securities and classes of assets have emerged in recent times and it is not always possible for an investor to specify them. In some cases, for instance, historical data of stocks are not available. In such cases, the uncertain returns of assets may be determined as interval numbers by using experts' knowledge.

In this paper, we propose an MVO2-like interval semi-absolute deviation model for portfolio selection, where the expected returns of securities are treated as interval numbers. Based on the concept of satisfaction index of interval inequality relation, we convert the interval semi-absolute deviation portfolio selection problem into two parametric linear programming problems.

This paper is organized as follows. In Section 2.1, we give some notations for interval numbers and briefly introduce some interval arithmetics. An order of relations over intervals is introduced. The concepts of satisfaction degree of interval inequality relations are given. Based on this concept, an approach to compare interval numbers is proposed. In Section 2.2, an approach is presented for estimating the intervals of rates of returns of securities. In Section 2.3, an interval absolute deviation model for portfolio selection is proposed. According to the approach proposed in Section 2.1, which concerns about comparing interval numbers, the interval portfolio selection problem is converted into a parametric linear programming problem with two parameters. In Section 3, an example is given to illustrate our approach. A few concluding remarks are finally given in Section 4.

## 2. Linear Programming Model with Interval Coefficients

### 2.1. Interval number and interval inequality

**Definition 2.1([2])**: Let $\circ \in \{+, -, \times, \div\}$ be a binary operation on $R$. If $a$ and $b$ are two closed intervals, then
$$a \circ b = \{x \circ y : x \in a, y \in b\}$$

defines a binary operation on the set of all the closed intervals. In the case of division, it is always assumed that 0 is not in $b$.

The operations on intervals used in this paper are as follows: for any two interval numbers $a = [\underline{a}, \overline{a}]$ and $b = [\underline{b}, \overline{b}]$,

$$a + b = [\underline{a} + \underline{b}, \overline{a} + \overline{b}], \ a - b = [\underline{a} - \overline{b}, \overline{a} - \underline{b}],$$
$$a \pm k = [\underline{a} \pm k, \overline{a} \pm k],$$

$$ka = k[\underline{a},\overline{a}] = \begin{cases} [k\underline{a}, k\overline{a}] & k \geq 0 \\ [k\overline{a}, k\underline{a}], & k < 0 \end{cases},$$

where $k$ is a real number.

An interval number can be viewed as a special fuzzy number whose membership function takes value 1 over the interval, and 0 anywhere else. For an interval number $a = [\underline{a}, \overline{a}]$, the median $m(a)$ and width $w(a)$ is defined by $m(a) = (\overline{a} + \underline{a})/2$ and $w(a) = (\overline{a} - \underline{a})/2$, respectively. The three operations of intervals are equivalent to the operations of addition, subtraction and scalar multiplication of fuzzy numbers via the extension principle. Ishibuchi and Tanaka suggested an order relation between two intervals as follows [6].

**Definition 2.2**: If intervals $a = [\underline{a}, \overline{a}]$ and $b = [\underline{b}, \overline{b}]$ are two intervals, the order relation between $a$ and $b$ is defined as

$$a \leq b \text{ if and only if } \underline{a} \leq \underline{b} \text{ and } \overline{a} \leq \overline{b}, \qquad (2.1)$$
$$a < b \text{ if and only if } a \leq b \text{ and } a \neq b \qquad (2.2)$$

For describing the interval inequality relation in detail, the following three concepts were introduced in [3]:

**Definition 2.3**: For any two interval numbers $a = [\underline{a}, \overline{a}]$ and $b = [\underline{b}, \overline{b}]$, there is an interval inequality relation $a \prec b$ between the two interval numbers $a$ and $b$ if and only if $m(a) \leq m(b)$. Furthermore, if $\overline{a} \leq \underline{b}$, we say the interval inequality relation $a \prec b$ between $a$ and $b$ is optimistic satisfactory; if $\overline{a} > \underline{b}$, we say the interval inequality relation $a \prec b$ between $a$ and $b$ is pessimistic satisfactory.

**Definition 2.4**: For any two interval numbers $a = [\underline{a}, \overline{a}]$ and $b = [\underline{b}, \overline{b}]$, if the interval inequality relation between them is pessimistic satisfactory, the pessimistic satisfaction index of the interval inequality relation $a \prec b$ can be defined as

$$PSD(a \prec b) = 1 + \frac{\underline{b} - \overline{a}}{w(a) + w(b)}. \qquad (2.3)$$

**Definition 2.5**: For any two interval numbers $a = [\underline{a}, \overline{a}]$ and $b = [\underline{b}, \overline{b}]$, if the interval inequality relation between them is optimistic satisfactory, the optimistic satisfaction index of the interval inequality relation $a \prec b$ can be defined as

$$OSD(a \prec b) = \frac{\underline{b} - \overline{a}}{w(a) + w(b)}. \qquad (2.4)$$

**Remark 2.1:** It is easy to see that $PSD(a \prec b) \geq 0$ if and only if $m(a) \leq m(b)$, $PSD(a \prec b) = 0$ if and only if $m(a) = m(b)$, and $OSD(a \prec b) = 0$ if and only if $PSD(a \prec b) = 1$.

Since $\overline{b} > \underline{a}$ implies that there may be some possibility for $b$ to be greater than $a$, it can not be said that the definition 2.2 and 2.3 contain all possibilities for interval inequality to hold. Therefore, we introduce an inclusive concept of interval inequality relation and satisfaction index.

**Definition 2.6**: For any two interval numbers $a = [\underline{a}, \overline{a}]$ and $b = [\underline{b}, \overline{b}]$, there is an interval inequality relation $a \preceq b$ between the two interval numbers $a$ and $b$ if and only if $\bar{b} > \underline{a}$. The satisfaction index of the interval inequality relation $a \preceq b$ is defined as

$$SD(a \preceq b) = \max\left\{\frac{\bar{b} - \underline{a}}{\bar{a} - \underline{a} + \bar{b} - \underline{b}}, 0\right\}. \qquad (2.5)$$

From (2.5), it follows that for every $\alpha \in (0, \infty)$, $SD(a \preceq b) \geq \alpha$ if and only if $\bar{b} - \underline{a} \geq \alpha(\bar{a} - \underline{a} + \bar{b} - \underline{b})$, i.e., $(1-\alpha)\bar{b} + \alpha\underline{b} \geq (1-\alpha)\underline{a} + \alpha\bar{a}$, and $SD(a \preceq b) \geq 1$ if and only if $\underline{b} \geq \bar{a}$. But, $OSD(a \preceq b) \geq 0$ if and only if $\underline{b} \geq \bar{a}$.

**Remark 2.2**: The possibility degree of $a \preceq b$ in [7] was defined by

$$PD(a \preceq b) = \min\left\{\max\left\{\frac{\bar{b} - \underline{a}}{\bar{a} - \underline{a} + \bar{b} - \underline{b}}, 0\right\}, 1\right\}.$$

According to definitions of the pessimistic, the optimistic satisfaction indices and possibility degree, we can see that the range of the pessimistic satisfaction index and possibility degree can be [0, 1), and the range of the optimistic satisfaction index and the satisfaction index can be [0, ∞). Our proposed satisfaction index is more general concept than the pessimistic, optimistic satisfaction indices of [3] and possibility degree of [7]. This fact is illustrated by following example. Suppose that possible relations between interval numbers $a$ and $b$ are such as following figures:

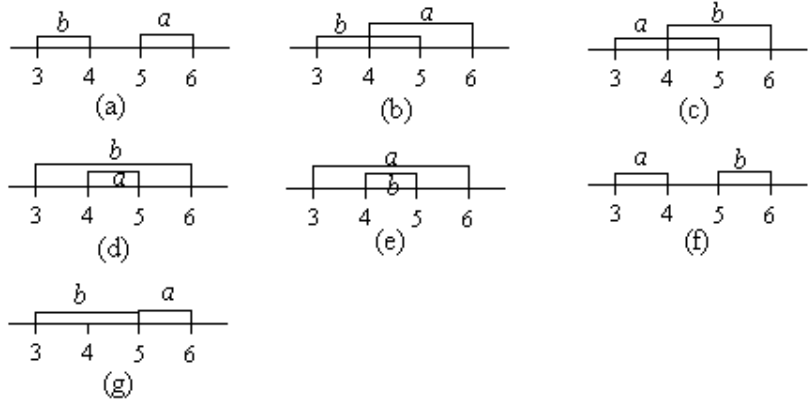

Then quadruple numbers consisting of corresponding optimistic and pessimistic satisfaction indices, possibility degree and satisfaction index are (a) (1, 1, 1, 1.5), (b) (0, 0.5, 0.75, 0.75), (c) (0, 0, 0.25, 0.25), (d) (0, 0, 0.5, 0.5), (e) (0, 0, 0.5, 0.5), (f) (0, 0, 0, 0) and (g) (0, 1, 1, 1), respectively. Obviously, our satisfaction index represents more precisely the interval inequality relations than others. The larger value of the satisfaction index is, the better is satisfied the interval inequality relation.

## 2.2 The Expected Return Intervals of Securities

It is well known that future returns of securities cannot be accurately predicted in any emerging securities market. Traditionally, researchers consider the arithmetic mean of historical returns as the expected return of the security. So the expected return of the security is a crisp value in this way. However for this technique, two main problems need to be solved [3]:

(1) If the time horizon of the historical data of a security is very long, the influence of the earlier historical data is the same as that of the recent data. However, recent data of a security most often indicate that the performance of a corporation is more important in recent data than in the earlier historical data.

(2) If the historical data of a security are not enough, one cannot accurately estimate the statistical parameters, due to data scarcity. Considering these two problems, perhaps it is a good idea to consider the expected return of a security as an interval number, rather than a crisp value, based on the arithmetic mean of historical data. Investors may make use of a corporation's financial reports and the security's historical data to determine the expected return interval's range. To determine the range of change in expected returns of securities, we will consider the following three factors as in [3]:

(1) Arithmetic mean: Although arithmetic means of returns of securities should not be expressed as expected returns directly, they are a good approximation. Denote the arithmetic mean return factor by $r_a$, which can be calculated with historical data.

(2) Historical return tendency: If recent returns of a security have been increasing, we can believe that the expected return of the security is greater than the arithmetic mean based on historical data. However, if recent returns of a security have been declining, we can assume that the expected return of the security is smaller than the arithmetic mean based on historical data. Denote the historical return tendency factor by $r_h$, which reflects the tendency of the return on the security. We can use the arithmetic mean of recent returns as $r_h$.

(3) Forecast of future returns of a security: The third factor influencing the expected return of a security is its estimated future returns. Denote the forecast return factor by $r_f$.

Estimation of $r_f$ requires some forecasts based on the financial reports and experts' individual experiences. Based on the above three factors, we can derive lower and upper limits of the expected return of the security. We can put the minimum of the three factors, $r_a$, $r_h$ and $r_f$, as the lower limit of the expected return, while we can put the maximum values of the three factors $r_a$, $r_h$ and $r_f$ as the upper limit of the expected return of the security.

## 2.3 The Interval Programming Models for Portfolio Selection

Assume that an investor wants to allocate his wealth among *n* risky assets offering random rates of returns and a risk-free asset offering a fixed rate of return. We introduce some notations as follows.

$\tilde{r}_j$ : the expected rate of return interval of risky asset $j$ ($j = 1, 2, \cdots, n$);

$r_{n+1}$ : the rate of return of risk-free asset $n + 1$;

$x_i$ : the proportion of the total investment devoted to risky asset $i$ ($i = 1, 2, \cdots, n$) or risk-free asset $n + 1$;

$x_i^0$ : the proportion of the risky asset $i$ ($i = 1, 2, \cdots, n$) or risk-free asset $n+1$ owned by the investor;

$r_{tj}$ : the historical rate of return of risky asset $j$ ($j = 1, 2, \cdots, n$), $t$ ($t = 1, 2, \cdots, T$);

$k_i$ : the rate of transaction costs for the asset $i$ ($i = 1, 2, \cdots, n + 1$);

$u_i$ : the upper bound of proportion of the total investment devoted to risky asset $i$ ($i = 1, 2, \cdots, n$) or risk-free asset $n + 1$.

We use a V shape function to express the transaction costs, so the transaction costs of the asset $i$ ($i = 1, 2, \cdots, n + 1$) can be denoted by

$$C_i(x_i) = k_i |x_i - x_i^0| \qquad (2.6)$$

So the total transaction costs of the portfolio $x = (x_1, x_2, \cdots, x_n, x_{n+1})$ can be denoted by

$$C(x) = \sum_{i=1}^{n+1} C_i(x_i) = \sum_{i=1}^{n+1} k_i |x_i - x_i^0| \qquad (2.7)$$

Denote

$$r_{aj} = \frac{1}{T} \sum_{t=1}^{T} r_{tj} \qquad (2.9)$$

The uncertain expected return of the risky asset $j$ ($j = 1, 2, \cdots, n$) can be represented as the following interval number:

$$\tilde{r}_j = [\underline{r}_j, \overline{r}_j] = [\min\{r_{aj}, r_{hj}, r_{fj}\}, \max\{r_{aj}, r_{hj}, r_{fj}\}], \qquad (2.10)$$

where $r_{aj}$ is the arithmetic mean factor of risky asset $j$, $r_{hj}$ is the historical return tendency factor of risky asset $j$ and $r_{fj}$ is the forecast return factor of risky asset $j$. They can be derived by using the above method. So the expected return interval of portfolio $x = (x_1, x_2, \cdots, x_{n+1})$ in the future can be represented as

$$\hat{r}(x) = \sum_{j=1}^{n} \tilde{r}_j x_j + r_{n+1} x_{n+1} \qquad (2.11)$$

After removing the transaction costs, the net expected return interval of portfolio $x = (x_1, x_2, \cdots, x_{n+1})$ can be represented as

$$\tilde{r}(x) = \sum_{j=1}^{n} \tilde{r}_j x_j + r_{n+1} x_{n+1} - \sum_{i=1}^{n+1} k_i |x_i - x_i^0| \qquad (2.12)$$

Because the expected returns on securities are considered as interval numbers, we may consider the semi-absolute deviation of the rates of return of portfolio $x$ below the expected return over all the past periods as an interval number too.

Since the expected return interval of portfolio $x = (x_1, x_2, \cdots, x_{n+1})$ is

$$\hat{r}(x) = [\sum_{j=1}^{n} \underline{r}_j x_j + \underline{r}_{n+1} x_{n+1}, \sum_{j=1}^{n} \overline{r}_j x_j + \underline{r}_{n+1} x_{n+1}] \qquad (2.13)$$

we can get the semi-absolute deviation interval of return of portfolio $x$ below the expected return over the past period $t$, $t = 1, 2, \cdots, T$. It can be represented as

$$\tilde{w}_t(x) = \left[ \max\left\{ \sum_{j=1}^{n} (\underline{r}_j - r_{tj}) x_j, 0 \right\}, \max\left\{ \sum_{j=1}^{n} (\overline{r}_j - r_{tj}) x_j, 0 \right\} \right]. \qquad (2.14)$$

Then the average value of the semi-absolute deviation interval of return of portfolio $x$ below the uncertain expected return over all the past periods, can be represented as

$$\tilde{w}(x) = \frac{1}{T} \sum_{t=1}^{T} \tilde{w}_t(x) \qquad (2.15)$$

We use $\tilde{w}(x)$ to measure the risk of portfolio $x$. Suppose that the investor wants to maximize the return of a portfolio after removing the transaction costs within some given level of risk. If the risk tolerance interval $\tilde{w} = [\underline{w}, \overline{w}]$ is given, the mathematical formulation of the portfolio selection problem is

(ILP)  
$$\max \ \tilde{r}(x) = \sum_{j=1}^{n} \tilde{r}_j x_j + r_{n+1} x_{n+1} - \sum_{i=1}^{n+1} k_i |x_i - x_i^0|$$
$$\text{s.t.} \ \tilde{w}(x) \preceq [\underline{w}, \overline{w}],$$
$$\sum_{i=1}^{n+1} x_i = 1, \ 0 \leq x_j \leq u_j, j = 1, 2, \cdots, n+1$$

where $\underline{w}$ represents the pessimistic tolerated risk level, and $\overline{w}$ represents the optimistic tolerated risk level.

(ILP) is an optimization problem with interval coefficients and, therefore, techniques of classical linear programming can not be applied unless the above interval optimization problem is reduced to a standard linear programming structure. In the following, we perform this conversion.

We introduce the order relation in the interval objective function of (ILP). Based on the concept of satisfaction index proposed by us in Section 2.1, the interval inequality relation $\tilde{w}(x) \preceq [\underline{w}, \overline{w}]$ in (ILP) is expressed by a crisp inequality. The crisp inequality equivalent to the interval constraint condition $\tilde{w}(x) \preceq [\underline{w}, \overline{w}]$ can be represented as follows:

$$SD(\tilde{w}(x) \preceq [\underline{w}, \overline{w}]) \geq \alpha \qquad (2.16)$$

Then the interval linear programming problem (ILP) can be represented by interval linear programming problem in which the objective function is interval number and the constraint conditions are crisp equality and inequalities. The interval objective function linear programming problem is represented as follows:

(IP) $\quad \max \ \tilde{r}(x) = \sum_{j=1}^{n} \tilde{r}_j x_j + r_{n+1} x_{n+1} - \sum_{i=1}^{n+1} k_i |x_i - x_i^0|,$

$\quad$ s.t. $SD(\tilde{w}(x) \preceq [\underline{w}, \overline{w}]) \geq \alpha,$

$$\sum_{i=1}^{n+1} x_i = 1, \ 0 \leq x_j \leq u_j, j = 1, 2, \cdots, n+1,$$

where satisfaction index $\alpha \in [0, \infty)$ is given by the investor.

Denote $F$ as the feasible set of (IP).

**Definition 2.7**: $x \in F$ is a satisfactory solution of (IP) if and only if there is no other $x' \in F$ such that $\tilde{r}(x) \prec \tilde{r}(x')$.

By Definition 2.7, the satisfactory solution of (IP) is equivalent to the non-inferior solution set of the following bi-objective programming problem for given satisfaction index $\alpha \in [0, \infty)$:

$\quad \max \ \bar{r}(x) = \sum_{j=1}^{n} \bar{r}_j x_j + r_{n+1} x_{n+1} - \sum_{i=1}^{n+1} k_i |x_i - x_i^0|$

(BLP) $\quad \max \ \underline{r}(x) = \sum_{j=1}^{n} \underline{r}_j x_j + r_{n+1} x_{n+1} - \sum_{i=1}^{n+1} k_i |x_i - x_i^0|$

$\quad$ s.t. $SD(\tilde{w}(x) \preceq [\underline{w}, \overline{w}]) \geq \alpha,$

$$\sum_{i=1}^{n+1} x_i = 1, \ 0 \leq x_j \leq u_j, j = 1, 2, \cdots, n+1$$

By the multi-objective programming theory, the non-inferior solution to (BLP) can be generated by solving the following parametric linear programming problem:

(PLP) $\quad \max \ \underline{r}(x) = \sum_{j=1}^{n} (\lambda \underline{r}_j + (1-\lambda) \bar{r}_j) x_j + r_{n+1} x_{n+1} - \sum_{i=1}^{n+1} k_i |x_i - x_i^0|$

$\quad$ s.t. $SD(\tilde{w}(x) \preceq [\underline{w}, \overline{w}]) \geq \alpha,$

$$\sum_{i=1}^{n+1} x_i = 1, \ 0 \leq x_j \leq u_j, j = 1, 2, \cdots, n+1$$

Introducing the concrete form of $SD(\tilde{w}(x) \preceq [\underline{w}, \overline{w}])$, (PLP) may be rewritten as follows:

(PLP1) $\max \ \underline{r}(x) = \sum_{j=1}^{n} (\lambda \underline{r}_j + (1-\lambda) \bar{r}_j) x_j + r_{n+1} x_{n+1} - \sum_{i=1}^{n+1} k_i |x_i - x_i^0|$

s.t. $(1-\alpha) \dfrac{1}{T} \sum_{t=1}^{T} \max \left\{ \sum_{j=1}^{n} (\underline{r}_j - r_{tj}) x_j, 0 \right\} + \alpha \dfrac{1}{T} \sum_{t=1}^{T} \max \left\{ \sum_{j=1}^{n} (\bar{r}_j - r_{tj}) x_j, 0 \right\} \leq (1-\alpha) \overline{w} + \alpha \underline{w},$

$$\sum_{i=1}^{n+1} x_i = 1, \ 0 \leq x_j \leq u_j, j = 1, 2, \cdots, n+1$$

To solve (PLP1), we consider the following transformation. First, we introduce a new variable $x_{n+2}$ such that

$$x_{n+2} \geq \sum_{i=1}^{n+1} k_i |x_i - x_i^0|. \tag{2.17}$$

Let

$$d_i^+ = \frac{|x_i - x_i^0| + (x_i - x_i^0)}{2}, \quad d_i^- = \frac{|x_i - x_i^0| - (x_i - x_i^0)}{2} \tag{2.18}$$

$$\underline{y}_t^+ = \frac{\left|\sum_{i=1}^n (\underline{r}_i - r_{ti})x_i\right| + \sum_{i=1}^n (\underline{r}_i - r_{ti})x_i}{2} \tag{2.19}$$

$$\bar{y}_t^+ = \frac{\left|\sum_{i=1}^n (\bar{r}_i - r_{ti})x_i\right| + \sum_{i=1}^n (\bar{r}_i - r_{ti})x_i}{2}. \tag{2.20}$$

Then, (PLP1) is equivalent to the following standard linear programming problem.

(PLP2) $\quad \max \ \underline{r}(x) = \sum_{j=1}^n (\lambda \underline{r}_j + (1-\lambda)\bar{r}_j)x_j + r_{n+1}x_{n+1} - x_{n+2}$

$\text{s.t.} \ \alpha \frac{1}{T}\sum_{t=1}^T \bar{y}_t^+ + (1-\alpha)\frac{1}{T}\sum_{t=1}^T \underline{y}_t^+ \leq (1-\alpha)\bar{w} + \alpha\underline{w},$

$\sum_{i=1}^{n+1} k_i (d_i^+ + d_i^-) \leq x_{n+2},$

$\underline{y}_t^+ - \sum_{i=1}^n (\underline{r}_i - r_{ti})x_i \geq 0, \ t = 1,...,T,$

$\bar{y}_t^+ - \sum_{i=1}^n (\bar{r}_i - r_{ti})x_i \geq 0, \ t = 1,...,T,$

$d_i^+ - d_i^- = x_i - x_i^0, \ i = 1,...,n+1,$

$d_i^+ \geq 0, d_i^- \geq 0, \ i = 1,...,n+1,$

$\underline{y}_t^+ \geq 0, \bar{y}_t^+ \geq 0, \ t = 1,...,T,$

$\sum_{i=1}^{n+1} x_i = 1, \ 0 \leq x_j \leq u_j, j = 1, 2, \cdots, n+1.$

One can use several algorithms of linear programming to solve (PLP2) efficiently, for example, the simplex method. So we can solve the original portfolio selection problem (ILP1) by solving (PLP2).

**Remark 2.3.** In [3], the interval inequality was replaced by equality $SD(\tilde{w}(x) \preceq [\underline{w}, \bar{w}]) = \alpha$,

while by inequality (2.16) in this paper. For fixed satisfaction index $\alpha$ and parameter $\lambda$, let $F(\alpha)$ and $V(\alpha,\lambda)$ denote the feasible set and optimal value of (PLP2), respectively. It is easy to see from the definition of satisfaction index that $\alpha_1 < \alpha_2$ implies $F(\alpha_2) \subset F(\alpha_1)$ and $V(\alpha_1,\lambda) \geq V(\alpha_2,\lambda)$. And, from the construction of objective function, it follows that $\lambda_1 < \lambda_2$ implies $V(\alpha,\lambda_1) \geq V(\alpha,\lambda_2)$. Therefore, $V(\alpha,\lambda)$ is non-increasing with $\alpha$ and $\lambda$. This means that the smaller satisfaction index is, the larger risk and, in turn, the larger return is. Given $\alpha \in [\underline{\alpha}, \overline{\alpha}]$ and $\lambda \in [0,1]$, we have $V(\overline{\alpha},1) \leq V(\alpha,\lambda) \leq V(\underline{\alpha},0)$, where $\underline{\alpha} > 0$ and $\underline{\alpha} < \overline{\alpha}$. For decision of appropriate $\alpha$ (or $\lambda$), we may use the grey comprehensive evaluation method [8] combining AHP (Analytic Hierarchal Process)[9] and TOPSIS (Technique for Ordered Preference by Similarity to Ideal Solution)[10] with attributes such as return and risk for given $\lambda$ (or $\alpha$).

## 3. Numerical Example

In this section, we suppose that an investor chooses 6 componential stocks and a risk-less asset for his investment. The rate of return of the risk-less asset is 0.0014 per month. We collected historical data of the 6 stocks during 8 periods, using one month as a period. Because the arithmetical methods do not produce good estimates of the actual returns that the investor will receive in the future, we forecasted $r_{fj}$, the return rate of risky asset $j$, according to Wavelet-Grey-SVR-Markov prediction method and obtain the expected rate of return interval of each stock. The historical return tendency $r_{hj}$ was obtained by $r_{hj} = \frac{1}{m}\sum_{t=T-m+1}^{T} r_{tj}$, where $m$ is amount of the most recent periods (we took $m=5$). The expected rate of return intervals are given in Table 3.1.

Suppose the investor stipulates risk level interval $\widetilde{w} = [0.015, 0.040]$; by the method proposed in the above section, we can solve the portfolio selection problem by solving (PLP2). For the given risk level interval $\widetilde{w}$, more satisfactory portfolios can be generated by varying the values of the parameters $\lambda$ and $\alpha$ in (PLP2).

**Table 3.1** The expected rates of returns intervals

|  | Stock1 | Stock2 | Stock3 | Stock4 | Stock5 | Stock6 |
| --- | --- | --- | --- | --- | --- | --- |
| Lower return | 0.0838 | 0.0562 | 0.0220 | 0.0600 | 0.0450 | 0.0488 |
| Upper return | 0.1000 | 0.0898 | 0.0513 | 0.0760 | 0.1040 | 0.0780 |

The return intervals, the risk intervals and the values of parameters of portfolios are listed in Table 3.2. The corresponding portfolios are listed in Table 3.3 and Table 3.4.

**Table 3.2.** The return intervals, the risk intervals and the values of parameters of portfolios

| $\alpha =1$ | $\lambda$ | Return interval | Risk interval |
|---|---|---|---|
| Portfolio 1 | 0 | [0.0150, 0.0363] | [0.0603, 0.0946] |
| Portfolio 2 | 0.12 | [0.0150, 0.0363] | [0.0603, 0.0946] |
| Portfolio 3 | 0.24 | [0.0150, 0.0330] | [0.0684, 0.0931] |
| Portfolio 4 | 0.36 | [0.0150, 0.0308] | [0.0708, 0.0919] |
| Portfolio 5 | 0.48 | [0.0150, 0.0306] | [0.0710, 0.0918] |
| Portfolio 6 | 0.60 | [0.0150, 0.0306] | [0.0710, 0.0918] |
| Portfolio 7 | 0.72 | [0.0150, 0.0304] | [0.0710, 0.0917] |
| Portfolio 8 | 0.84 | [0.0150, 0.0304] | [0.0710, 0.0917] |
| Portfolio 9 | 0.96 | [0.0150, 0.0304] | [0.0710, 0.0917] |

| $\alpha =0.5$ | $\lambda$ | Return interval | Risk interval |
|---|---|---|---|
| Portfolio 1 | 0 | [0.0203, 0.0347] | [0.0720, 0.0975] |
| Portfolio 2 | 0.12 | [0.0214, 0.0336] | [0.0747, 0.0974] |
| Portfolio 3 | 0.24 | [0.0227, 0.0323] | [0.0777, 0.0968] |
| Portfolio 4 | 0.36 | [0.0227, 0.0323] | [0.0777, 0.0968] |
| Portfolio 5 | 0.48 | [0.0229, 0.0321] | [0.0781, 0.0966] |
| Portfolio 6 | 0.60 | [0.0229, 0.0321] | [0.0782, 0.0964] |
| Portfolio 7 | 0.72 | [0.0229, 0.0321] | [0.0782, 0.0964] |
| Portfolio 8 | 0.84 | [0.0229, 0.0321] | [0.0782, 0.0964] |
| Portfolio 9 | 0.96 | [0.0229, 0.0321] | [0.0782, 0.0964] |

**Table 3.3.** The allocation of portfolio 1, 2, 3, 4, 5, 6, 7, 8, 9 for $\alpha =1$

| portfolio | 1 | 2 | 3 | 4 | 5 | 6 | 7 | 8 | 9 |
|---|---|---|---|---|---|---|---|---|---|
| $\lambda$ | 0 | 0.12 | 0.24 | 0.36 | 0.48 | 0.60 | 0.72 | 0.84 | 0.96 |
| Stock 1 | 0.2378 | 0.2378 | 0.4462 | 0.5550 | 0.5639 | 0.5639 | 0.5716 | 0.5716 | 0.5716 |
| Stock 2 | 0.4317 | 0.4317 | 0.3424 | 0.2433 | 0.2342 | 0.2342 | 0.2235 | 0.2235 | 0.2235 |
| Stock 3 | 0 | 0 | 0 | 0.0364 | 0.0402 | 0.0402 | 0.0451 | 0.0451 | 0.0451 |
| Stock 4 | 0.0831 | 0.0831 | 0.1518 | 0.1608 | 0.1617 | 0.1617 | 0.1598 | 0.1598 | 0.1598 |
| Stock 5 | 0.2474 | 0.2474 | 0.0597 | 0.0046 | 0 | 0 | 0 | 0 | 0 |
| Stock 6 | 0 | 0 | 0 | 0 | 0 | 0 | 0 | 0 | 0 |
| Stock 7 | 0 | 0 | 0 | 0 | 0 | 0 | 0 | 0 | 0 |

**Table 3.4.** The allocation of portfolio 1, 2, 3, 4, 5, 6, 7, 8, 9 for $\alpha$ =0.5

| portfolio | 1 | 2 | 3 | 4 | 5 | 6 | 7 | 8 | 9 |
|---|---|---|---|---|---|---|---|---|---|
| $\lambda$ | 0 | 0.12 | 0.24 | 0.36 | 0.48 | 0.60 | 0.72 | 0.84 | 0.96 |
| Stock 1 | 0.6135 | 0.6877 | 0.8214 | 0.8214 | 0.8586 | 0.8733 | 0.8733 | 0.8733 | 0.8733 |
| Stock 2 | 0.2850 | 0.2692 | 0.1387 | 0.1387 | 0.0894 | 0.0671 | 0.0671 | 0.0671 | 0.0671 |
| Stock 3 | 0 | 0 | 0.0328 | 0.0328 | 0.0520 | 0.0596 | 0.0596 | 0.0596 | 0.0596 |
| Stock 4 | 0 | 0 | 0.0071 | 0.0071 | 0 | 0 | 0 | 0 | 0 |
| Stock 5 | 0.1014 | 0.0431 | 0 | 0 | 0 | 0 | 0 | 0 | 0 |
| Stock 6 | 0 | 0 | 0 | 0 | 0 | 0 | 0 | 0 | 0 |
| Stock 7 | 0 | 0 | 0 | 0 | 0 | 0 | 0 | 0 | 0 |

The investor may choose his own investment strategy from the portfolios according to his attitude towards the securities' expected returns and the degree of portfolio risk with which he is comfortable. If the investor is not satisfied with any of these portfolios, he may obtain more by solving the parametric linear programming problems (PLP2) for other values of parameter $\lambda$ and $\alpha$.

## 4. Conclusion

In [3], an approach was presented for estimating intervals of rates of returns of securities and the semi-absolute deviation risk function was extended to an interval case. They proposed an interval semi-absolute deviation model with no short selling and no stock borrowing in a frictional market for portfolio selection. In this paper, by introducing a concept of inclusive satisfaction index of the interval inequality relation, an approach to compare interval numbers is given. By using the approach, the interval semi-absolute deviation model can be converted into a parametric linear programming problem with two parameters. We represented the interval inequality by the satisfaction index inequality unlike equality of [3]. One can find a satisfactory solution to the original problem by solving the corresponding parametric linear programming problem. An investor may choose a satisfactory investment strategy according to an optimistic or pessimistic attitude by choosing proper values of parameter $\alpha$ and $\lambda$. The model can help the investor to find an efficient portfolio according to his/her preference.